\documentclass[pdflatex,sn-mathphys-num]{sn-jnl}
\usepackage{amsthm, amsmath, amssymb,latexsym}

\usepackage{graphicx}%
\usepackage{multirow}%
\usepackage{amsmath,amssymb,amsfonts}%
\usepackage{amsthm}%
\usepackage{mathrsfs}%
\usepackage[title]{appendix}%
\usepackage{xcolor}%
\usepackage{textcomp}%
\usepackage{manyfoot}%
\usepackage{booktabs}%
\usepackage{algorithm}%
\usepackage{algorithmicx}%
\usepackage{algpseudocode}%
\usepackage{listings}%

\newtheorem{theorem}{Theorem}[section]

\newtheorem{corollary}[theorem]{Corollary}
\newtheorem{proposition}[theorem]{Proposition}
\newtheorem{remark}[theorem]{Remark}
\newtheorem{definition}[theorem]{Definition}

 \newenvironment{dedication}
     {\vspace{0ex}\begin{quotation}\begin{center}\begin{em}}
     {\par\end{em}\end{center}\end{quotation}}
\newcommand{\nc}{\newcommand} 
\nc{\cH}{{\mathcal H}}
\nc{\cA}{{\mathcal A}}
\nc{\cG}{{\mathcal G}}
\nc{\cC}{{\mathcal C}}
\nc{\cO}{{\mathcal O}}
\nc{\cI}{{\mathcal I}}
\nc{\cB}{{\mathcal B}}
\nc{\cY}{{\mathcal Y}}
\nc{\cK}{{\mathcal K}} 
\nc{\cX}{{\mathcal X}}
\nc{\cS}{{\mathcal S}}
\nc{\cE}{{\mathcal E}}
\nc{\cF}{{\mathcal F}}
\nc{\cZ}{{\mathcal Z}}
\nc{\cQ}{{\mathcal Q}}
\nc{\cN}{{\mathcal N}}
\nc{\cP}{{\mathcal P}}
\nc{\cL}{{\mathcal L}}
\nc{\cM}{{\mathcal M}}
\nc{\cT}{{\mathcal T}}
\nc{\cW}{{\mathcal W}}
\nc{\cU}{{\mathcal U}}
\nc{\cJ}{{\mathcal J}}
\nc{\cV}{{\mathcal V}}
\nc{\bH}{{\mathbb H}}
\nc{\bA}{{\mathbb A}}
\nc{\bG}{{\mathbb G}}
\nc{\bC}{{\mathbb C}}
\nc{\bO}{{\mathbb O}}
\nc{\bI}{{\mathbb I}}
\nc{\bB}{{\mathbb B}}
\nc{\bY}{{\mathbb Y}}
\nc{\bK}{{\mathbb K}} 
\nc{\bX}{{\mathbb X}}
\nc{\bS}{{\mathbb S}}
\nc{\bE}{{\mathbb E}}
\nc{\bF}{{\mathbb F}}
\nc{\bZ}{{\mathbb Z}}
\nc{\bQ}{{\mathbb Q}}
\nc{\bN}{{\mathbb N}}
\nc{\bP}{{\mathbb P}}
\nc{\bL}{{\mathbb L}}
\nc{\bM}{{\mathbb M}}
\nc{\bT}{{\mathbb T}}
\nc{\bW}{{\mathbb W}}
\nc{\bU}{{\mathbb U}}
\nc{\bD}{{\mathbb D}}
\nc{\bJ}{{\mathbb J}}
\nc{\bV}{{\mathbb V}}
\nc{\bbZ}{{\mathbb Z}}
\nc{\bR}{{\mathbb R}}
\nc{\fr}{{\rightarrow}}
\nc{\co}{{\nabla}}
\newcommand{\la}{\longrightarrow}
\nc{\cu}{{\barline{\nabla}}}
\nc{\Aut}{\hbox{Aut}}
\newcommand{\sff}{\mathrm{I\!I}}

\begin{document}

\title{Projective structures and Hodge theory}

\author[1]{\fnm{Andrea} \sur{Causin}}\email{acausin@uniss.it}
\author*[2]{\fnm{Gian Pietro} \sur{Pirola}}\email{gianpietro.pirola@unipv.it}

\affil[1]{\orgdiv{DADU}, \orgname{Universit\`a di Sassari}, \orgaddress{\street{piazza Duomo 6}, \city{Alghero SS}, \postcode{07041}, \country{Italy}}}

\affil*[2]{\orgdiv{Dipartimento di Matematica}, \orgname{Universit\`a di Pavia}, \orgaddress{\street{via Ferrata 5}, \city{Pavia PV}, \postcode{27100},  \country{Italy}}}


\abstract{Every compact Riemann surface $X$ admits a natural projective structure $p_u$ as a consequence of the uniformization theorem. In this work we describe the construction of another natural projective structure on $X$, namely the Hodge projective structure $p_h$, related to the second fundamental form of the period map. We then describe how projective structures correspond to $(1,1)$-differential forms on the moduli space of projective curves and, from this correspondence, we deduce that $p_u$ and $p_h$ are not the same structure.}



\maketitle

  \begin{dedication}
{This article is dedicated to} {Alberto} {Collino}
\end{dedication}

\section{Introduction} 
A projective structure on a Riemann surface is a maximal atlas 
with linear fractional functions (M\"{o}bius transformations) as coordinate changes.
As a consequence of the uniformization theorem there exists a preferred structure
$p_u$ on every Riemann surface. 
The space of projective structures on a compact Riemann surface $X$ is a complex affine space $P_g$ 
having the space of holomorphic quadratic differentials on $X$ as corresponding vector space.
The structure $p_u$ can then be interpreted as a $\mathcal C^\infty$-section of an affine bundle
 $\cP$ over $\cM_g$, the moduli space of compact Riemann surfaces of genus $g$. 
 
The problem of comparing projective structures defined on families of compact Riemann surfaces 
is then addressed by passing to some related differential forms.
In particular, with this method we show that the canonical structure $p_u$ is different from
 the Hodge canonical structure, of which we sketch the construction.

 A second canonical projective (Hodge-)structure  $p_h$ on the compact Riemann surface $X$ is 
 indirectely constructed in \cite{CFG}. 
 In this paper, the second fundamental form corresponding to the period map 
 is represented by a meromorphic $2$-form on the product $X\times X$, 
 with poles only in the diagonal;
 the (Hodge) projective structure $p_h$ is constructed  from this $2$-form through the solutions of the schwarzian equation.
 In terms of differential equations, $p_h$ is associated to the Laplacian (commutative Hodge theory) while 
 $p_u$ is associated to a nonlinear differential equation (noncommutative Hodge theory).
 We then have two canonical sections $p_u$ and $p_h$ of the bundle $\cP$ 
 and the comparison between them is possible thanks to the interpretation 
 of the derivatives of the sections of $\cP$ as $(1,1)$-differential forms on the moduli space $\cM_g$.
 It is a classical result that the form $w_u$ corresponding to $p_u$ is the K\"{a}hler form 
 related to the Weil-Petersson metric. In \cite{BCFP} it is proven that the form $w_h$ corresponding to $p_h$ is the pull-back of the K\"{a}hler form  to the Siegel space through the period map. 
 In a subsequent paper \cite{FPT} the connection between projective structures and differential forms 
 is specified; furthermore, this relation is used in \cite{BFPT} in order to clarify the relation between the form $w_h$ and the section $p_h$ by means of the Quillen metrics.
 The geometry of the atlas corresponding to the projective Hodge structure is still completely mysterious.
 
 \bigskip
 
This paper is based on the plenary conference {\em ``Strutture proiettive e teoria di Hodge''} given by the second author at the XXII Congress of the {\em Unione Matematica Italiana}, held in Pisa in  September 5, 2023.


\section{Cartography and Global Positioning System (GPS)}

\subsection{Charts}
The fundamental problem of cartography is the representation of the Earth's surface at a small scale
 by means of charts, and the second problem is the {\em change of charts}, 
 that is the difficulty of producing easily usable coordinate map transitions. 
 
\begin{figure}[h!]
\centerline{\includegraphics[width=0.9\textwidth]{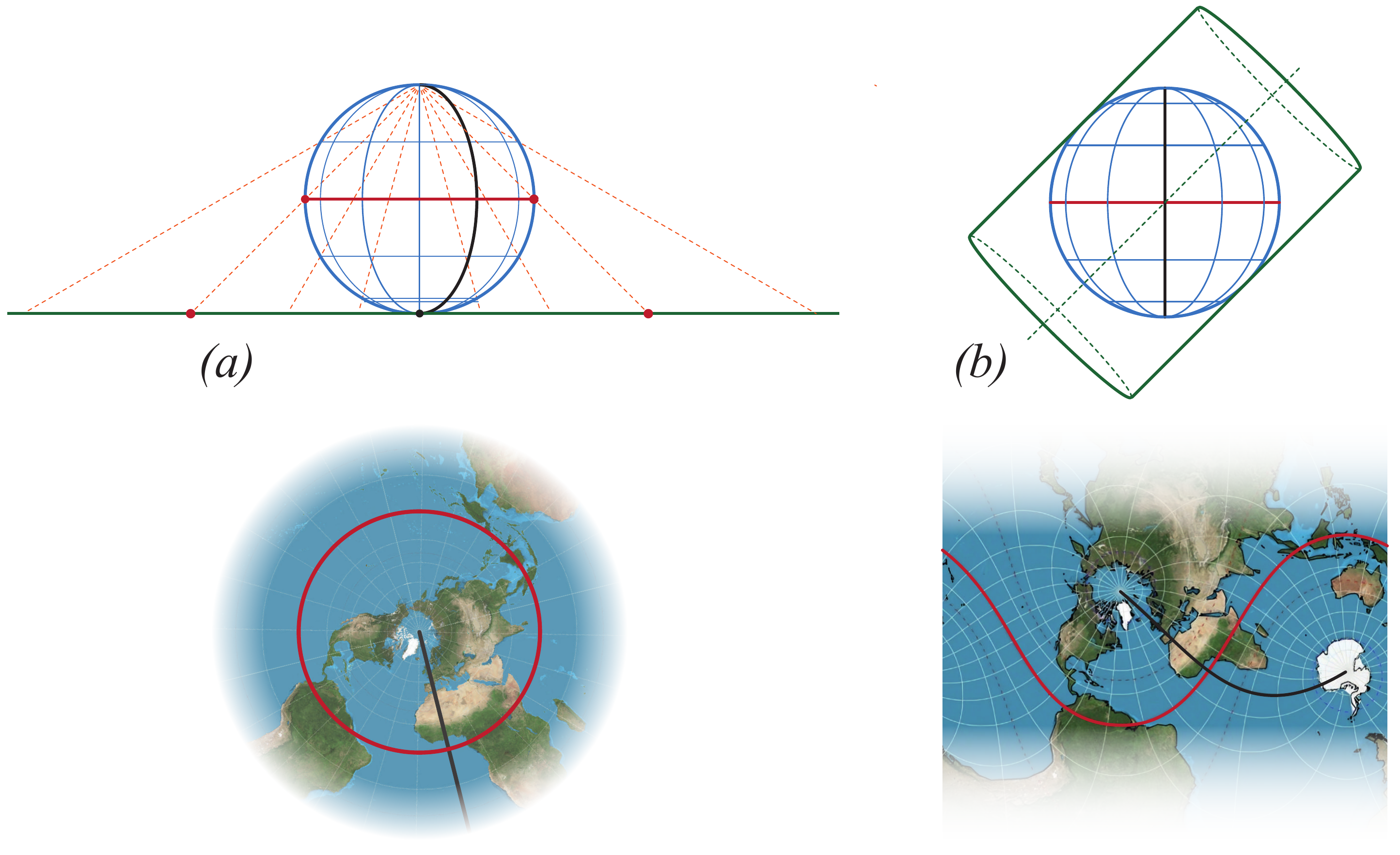}}
\caption{(a) Stereographic projection from a pole. (b) Oblique Mercator projection. }
\label{fig1}
\end{figure} 
 
Stereographic projections from the sphere to the plane or to the cylinder 
were classically used (Figure \ref{fig1}, a). The most famous charts are perhaps the ones
constructed in $1596$ by the Flemish geographer and cartographer Gerard Mercator
via conformal cylindrical cartographic projections (Figure \ref{fig1}, b). 
A curious remark is that Mercator maps are related to the complex exponential (see {\em e.~g.} \cite{P}).

An attempt to go beyond the stereographic projections
is made by Lagrange in his two memoirs of $1779$ and $1781$ both with title {\em ``Sur la construction de cartes g\'eographiques''} \cite{L}.

 \subsection{Atlases}
A non-secondary aspect and an important geometric and topological problem on every differentiable manifold 
is that of finding atlases with simple changes of coordinates. We will treat here only the case of surfaces and we recall 
that the general existence of conformal atlases (in the analytic case) can be traced back to Gauss and this has surely been 
 a source of inspiration of Riemann's work. 
 
 \begin{figure}[h!]
\centerline{\includegraphics[width=0.7\textwidth]{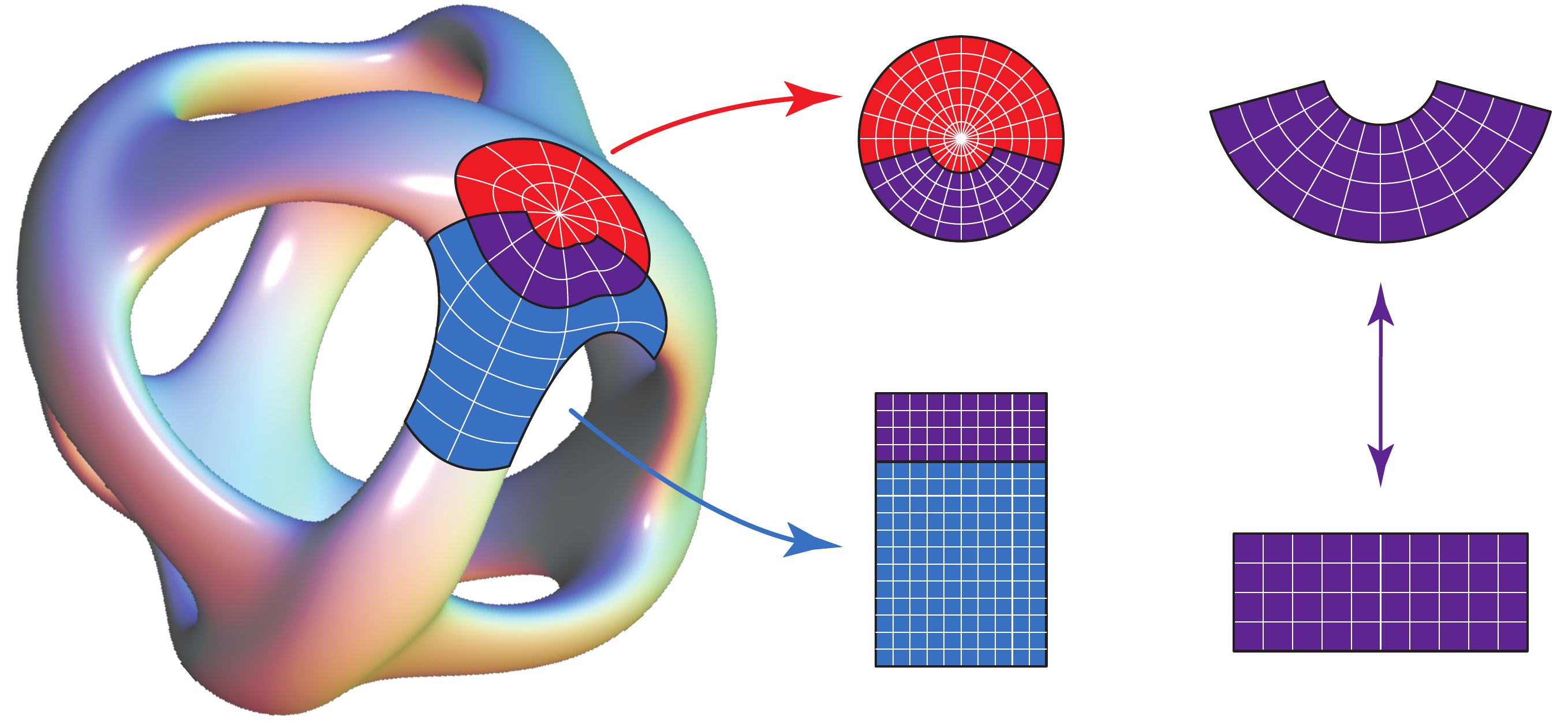}}
\caption{Two charts on a (Riemann) surface }
\label{fig2}
\end{figure} 
 
 We now call {\em Riemann surface} every connected orientable surface $X$ 
 with holomorphic ($=$ conformal) changes of charts (Figure \ref{fig2}). A natural question is to 
 ask if there exist simple changes of charts.

\begin{figure}[h!]
\centerline{\includegraphics[width=0.6\textwidth]{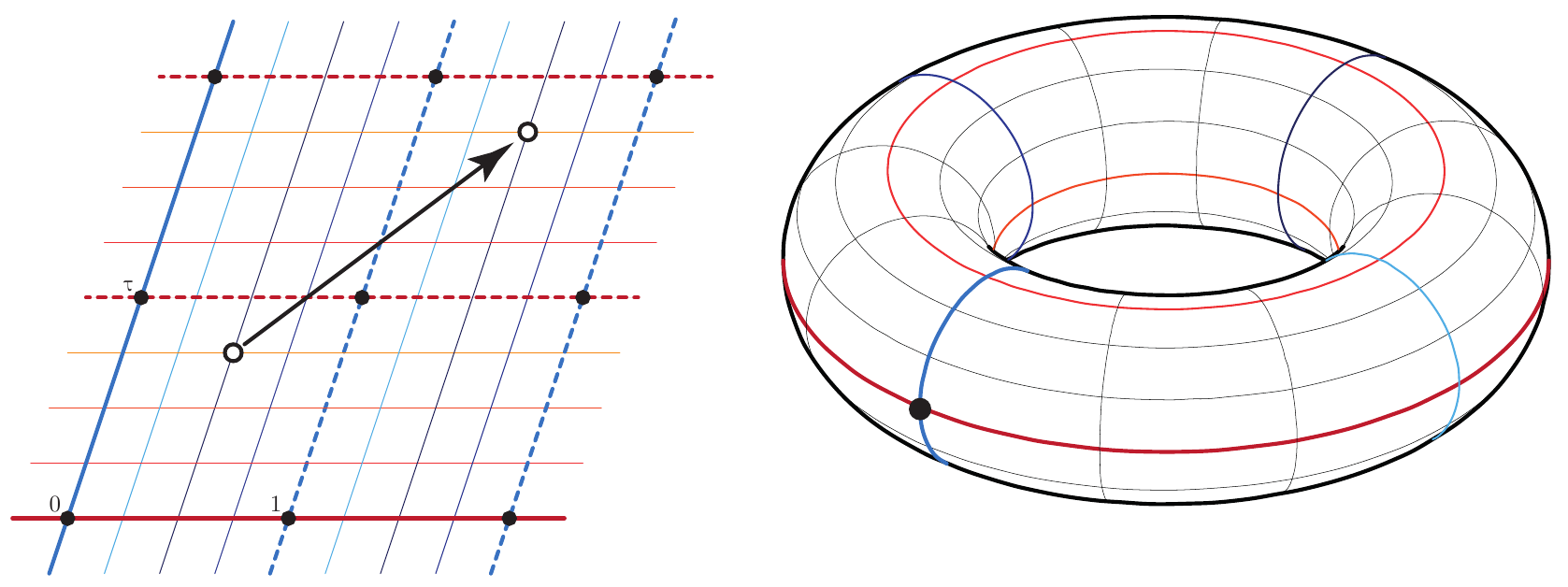}}
\caption{Affine transformations produce tori. }
\label{fig3}
\end{figure}

 The simplest coordinate changes are the affine transformations
$$z\mapsto az+b$$
but they prove inadequate, as they impose severe topological restrictions.
The existence of an affine atlas on a compact surface $X$, for instance, implies that $X$ is a torus (Figure \ref{fig3})
$$X\equiv S^1\times S^1,$$ where $S^1$ is the circumference.

\subsection{Uniformization}

The uniformization theorem solves the problem of finding simple atlases for surfaces.

\begin{theorem}[Riemann \cite{R}, Poincar\'e \cite{Poincare}, Koebe \cite{K} and see \cite{H}] Every simply connected Riemann surface is biholomorphic  to either:
\begin{enumerate}
\item the Riemann sphere $\bS^2=\bP^1=\bC\cup\{\infty\}$, 
\item the complex plane $\bC,$ 
\item the Poincar\'e half-plane $\bH=\{z\in \bC \ : \ \hbox{Im}(z)>0\}$ 

equivalent to the disk $\bD=\{z\in \bC \ : \ |z|<1\}.$
 \end{enumerate}\end{theorem}

The universal cover $$u: U\to X$$ of a Riemann surface $X$ is a parametrization of $X$ {\em from above}
 and can be interpreted as a Global Positioning System (GPS). 

\begin{figure}[h!]
\centerline{\includegraphics[width=0.9\textwidth]{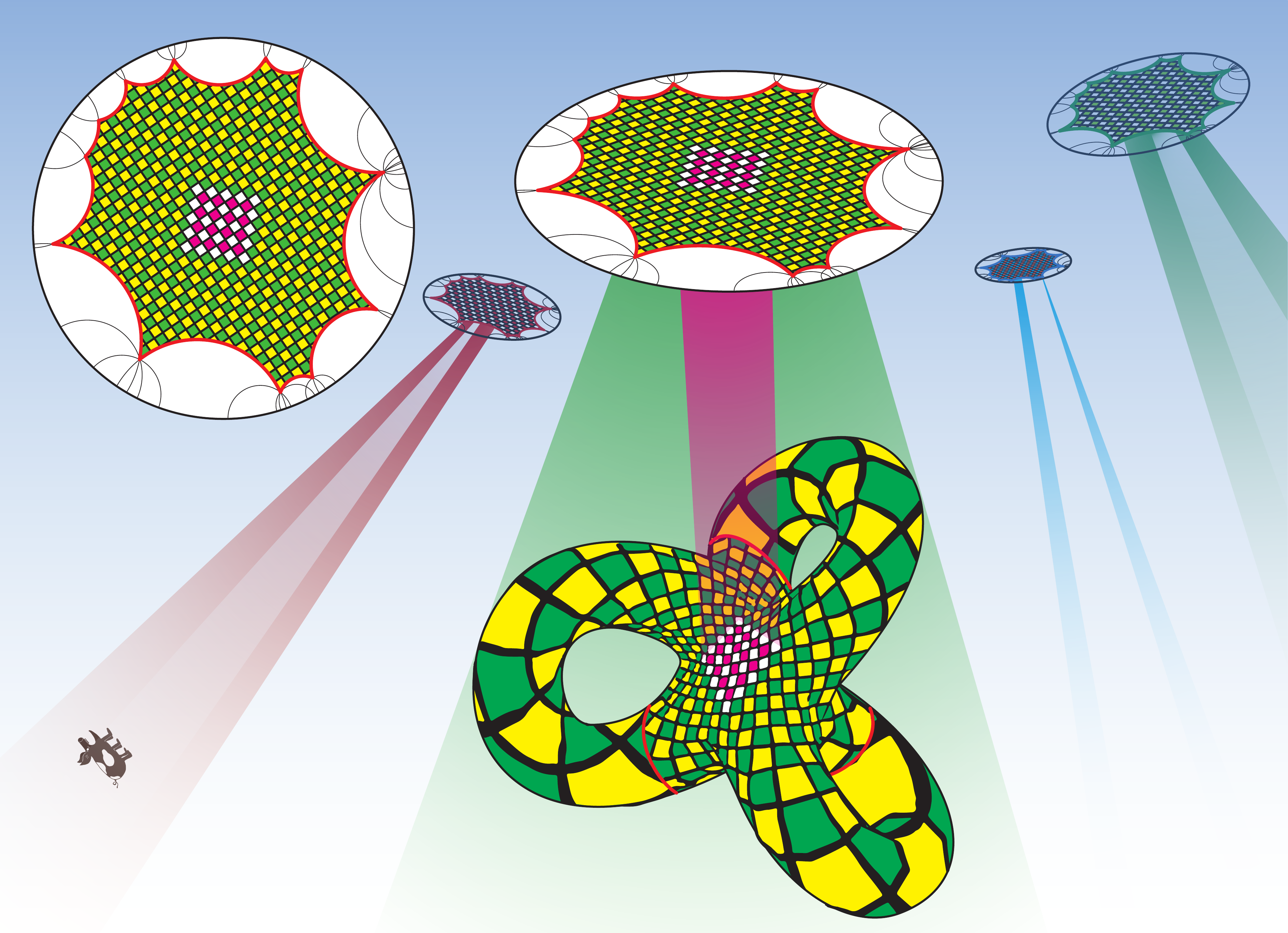}}
\caption{The universal cover parametrizes the surface from above. }
\label{fig4}
\end{figure}

\subsection{Automorphisms}

The group of biholomorphisms of $\bP^1$, known as the M\"{o}bius group, 
 is the group of linear fractional transformations 
$$\Aut(\bP^1)=\left\{ g_M(z)=\frac{az+b}{cz+d} \ : \  M=\begin{pmatrix} a&b\\ c&d\end{pmatrix}\in GL(2,\bC) \ ,\ \det(M)=1\right\}$$
where we interpret $\bP^1$ as $\bC \cup\{\infty\}.$ 

\noindent Since $g_M=g_N$ if and only if $ M=\pm N$, we have 
$$Aut(\bP^1)\cong \frac{SL(2,\bC)}{\{\pm I\}}= \bP GL(2,\bC).$$
The automorphisms of $\bC$ are all the affine transformations
 $$Aut(\bC)=\{f(z)= az+b \ : \ a,b \in\bC, \ a\neq 0\},$$ 
 and they can be identified with the linear fractional transformations 
 fixing the point at infinity.

\noindent Lastly, the automorphisms of $\bH$ are the restrictions to the haf-plane of the 
 real linear fractional transformations 
 $$ Aut (\bH)=\bP GL(2,\bR).$$
In conclusion,
\begin{enumerate}
\item $Aut(\bP^1)\cong \bP GL(2,\bC)$
\item $Aut(\bC) \subset \bP GL(2,\bC)$
\item $Aut(\bH)=\bP GL(2,\bR)\subset \bP GL(2,\bC).$
\end{enumerate}

\section{Projective structures}
We would like to begin this section by quoting E. Cartan \cite{Cartan} 
who explains in a way that could  not be clearer:
 \bigskip
 
{\em Une vari\'et\'e (ou espace) \`a connexion projective est une 
vari\'et\'e num\'erique qui, au voisinage imm\'ediat de chaque point, 
pr\'esente tous les caract\`eres d'un espace projectif et dou\'ee de plus 
d'une loi permettant de raccorder en un seul espace projectif les deux petits morceaux  [...]}.

\subsection{Projective structures on surfaces}

A projective structure on a Riemann surface $X$ is a (maximal) atlas
 with linear fractional coordinate map transitions, 
 compatible with the conformal structure of $X$. More precisely: 

\begin{definition} A projective structure on $X$ is a maximal atlas 
$\{(U_\alpha, \phi_\alpha)\}_\alpha$ where $U_\alpha$ is an open subset of $X$, 
$$\phi_\alpha: U_\alpha \rightarrow \bP^1$$
is an injective holomorphic function, and 
$$\phi_\alpha\circ \phi_\beta^{-1}: \phi_\beta(U_\alpha\cap U_\beta)\to \phi_\alpha(U_\alpha\cap U_\beta):$$
is a restriction of an element
$$g_{\alpha\beta}\in \bP GL(2,\bC).$$
 \end{definition}

The local sections of the universal cover $u: U\rightarrow X$ 
provide a natural atlas with coordinate map changes in $Aut(U)$
and the uniformization theorem implies that $Aut(U)\subset \bP GL(2,\bC)$.

Therefore:

\begin{proposition} Every Riemann surface admits a projective structure. \end{proposition}

\subsection{Projective structures on compact surfaces}

From now on, we will assume that $X$ is a compact Riemann surface of genus $g\geq 2$.
The uniformization theorem provides a natural projective structure $p_u(X)$ on $X$; 
this structure is, however, far from being unique:

\begin{theorem}[Gunning, \cite{Gunning}] \label{strutt} Let $X$ be a compact Riemann surface of genus $g\geq 2$. 

\noindent The space $P_X$ of  projective structures on $X$ is a complex affine space 
corresponding to the vector space $H^0(X, 2K_X)$ of quadratic differentials on $X$. In particular
$$\dim_\bC P_X=3g-3.$$ 
\end{theorem}

\bigskip

The choice of $p_u(X)$ as origin for $P_X$ gives an isomorphism $P_X\cong H^0(X, 2K_X)$,
 that is  a vector space structure on $P_X$. 
 
 It should be remarked that even if $X_t$ is a holomorphic family of compact Riemann surfaces,
the corresponding family $p_{u,t}$ of complex structures on $X_t$ is, in general, 
 smooth ($\mathcal C^\infty$) but not holomorphic with respect to $t$.

\bigskip

According to Riemann, let $\mathcal M_g$ denote the (Riemann) moduli space 
of the conformal structures on $X$; 
we have that $\mathcal M _g$ is a (singular) algebraic variety 
of dimension $3g-3$. We will however neglect the singularities 
as we have several technical tools for treating them 
(Stack-Orbifold, level structures, or just the restriction to a smooth open subset).

Let us indicate by 
$$\varphi:T^\ast_{\mathcal M_g}\rightarrow \mathcal M_g$$ 
the cotangent bundle (as said, by neglecting singularities); 
its fiber over $[X]$ corresponds to $H^0(X, 2K_X).$
The space of projective structures 
$$\psi:\cP\to \cM_g$$ 
is an affine space of complex dimension $6g-6$ and with $T^\ast_{\cM_g}$ 
as corresponding vector bundle.
An important result to which we will refer later is that the fibration $\psi:\cP\to \cM_g$
does not admit holomorphic sections so, in particular, this says again that the section $p_u$ 
is not holomorphic.

 \bigskip

We now briefly discuss the theorem \ref{strutt}.
The linear fractional transformations 
$$f(z)= \frac{az+b}{cz+d}$$ 
can be characterized as the functions that nullify the Schwarzian derivative

$$ S(f)= \frac{f'''}{f'}-\frac{3}{2}\left(\frac{f'' }{f'}\right) ^2.$$

  \bigskip

We can also adopt a more intrinsic point of view: 
if $(X,\sigma)$ is a Riemann surface with a fixed projective structure and
$f:X\rightarrow \bC$ is holomorphic, then 

 $$ S(f)_\sigma= \left(\frac{f'''}{f'}-\frac{3}{2}\left(\frac{f'' }{f'}\right)^2 \right)dz^2$$ 
 
 is a meromorphic quadratic differential, well-defined on $X$.
 
 Given two projective structures $\sigma$ and $\sigma^\prime$ on $X$ and
 
  $${\rm id}: (X,\sigma)\to (X,\sigma')$$ then
 \begin{equation}\label{gls}  \sigma-\sigma' :=S({\rm id})_\sigma=\alpha  \end{equation}
  is a quadratic differential. 
  Moreover, given $\{\sigma,\alpha\} $, there exists $\sigma^\prime$ solving (\ref{gls}).

 \begin{remark}
 The Schwarzian differential operator is already present in the Lagrange memoirs 
 ``Sur la construction des cartes g\'eographiques'' we previously cited. \end{remark}

\section{ The Hodge projective structure}

A new natural projective structure on every compact Riemann surface of genus $g\geq 2$ 
was found by Elisabetta Colombo, Paola Frediani and Alessandro Ghigi 
by means of Hodge theory. It is a winding route, 
all the same, we will try to outline a synthetic path along it.

\bigskip
\subsection{The period matrix}

The Hodge theory on a Riemann surface $X$ is, in simple terms, the theory of 
the Laplacian on a surface.

As before, let $X$ be a compact Riemann surface of genus $g$.
 By associating a compatible metric, we have that the 
 Dirichlet form is conformally invariant: 
on a differential $1$-form $\omega\equiv pdx+q dy$
depending locally on the conformal coordinate $z=x+iy$,
 the Hodge star operator $\star$ is defined as
$$ \star \omega\equiv -qdx+pdy; $$ 

it is conformally invariant and allows to define the product ({\em Dirichlet form})

\begin{equation}\label{star} 
(\omega,\omega^\prime) =\int_X \omega\wedge \star \omega^\prime 
\end{equation}

on the space $\cA^1(X)$ of all differential $1$-forms on X.

Let $d$ denote the exterior differential; the real cohomology group defined as

$$ H^1(X,\bR) =\frac{\cA^1(X)_{closed}}{\cA^1(X)_{exact}}=
\frac {\{\omega\in \cA^1(X) 	\ : \ d\omega=0\}}{\{\omega\in \cA^1(X) \ : \  \omega=df, \ f\in \cC^\infty(X)\}}$$

captures the {\em potential theory} of $X$.

Fixing a class $\alpha\in H^1(X,\bR)$ we can employ 
the scalar product induced by (\ref{star}) in order to minimize 
the Dirichlet norm
$$ \cF(\omega)= \int_X \omega\wedge \star \omega $$
obtaining
$$ \omega_\alpha = \hbox{argmin} \left\{ \cF(\omega) \ : \ \omega \in\cA^1(X) , \  [ \omega ]=\alpha \right\}.  $$ 

It can be shown that for any $ \alpha\in H^1(X,\bR) $ the minimizer $\omega_\alpha$ exists, is unique and harmonic; 
more precisely, it is the real part 
$$\omega_\alpha={\rm Re}(\Omega_\alpha)$$ 
of a holomorphic 1-form 
$$ \Omega_\alpha \equiv f_\alpha(z)\, dz, \quad \hbox{ with } f_\alpha(z) \hbox{ holomorphic.}$$
This minimization procedure allows to select a precise subspace of holomorphic $1$-forms $H\subset{H^0(X, \Omega^1_X)}$ 
representing the cohomology $H^1(X,\bR)$.
\bigskip
 
\begin{figure}[h!]
\centerline{\includegraphics[width=0.4\textwidth]{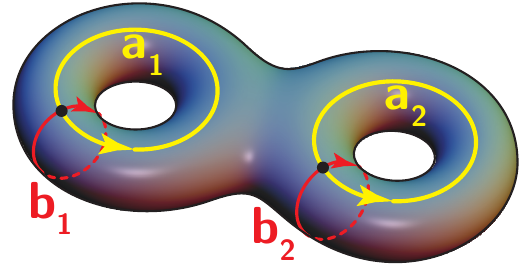}}
\caption{A symplectic basis for the (homology of the) genus $2$ surface.}
\label{fig5}
\end{figure}

 These minimizers can be organized as follows. 
 Fix a symplectic (with respect to the intersection coupling) basis 
$a_1,\dots , a_g , b_1,\dots , b_g$ of $H_1(X,\bZ) $ (Figure \ref{fig5})
and choose a basis of $H$

$$ \Omega_1, \dots , \Omega_g  $$

in such a way that 

$$\int_{a_i}\Omega_j = \delta_{i,j}= \left\{ \begin{array}{lr} 1 & \hbox{ if } i=j \\ 0 & \hbox{ if } i\neq j \end{array}\right. .$$
 
Then, still following Riemann, the period matrix is defined as  
$$\Omega(X)=(p_{i,j})=\left(\int_{b_i}\Omega_j\right).$$

 The matrix $\Omega(X)$ is symmetric and with positive definite imaginary part, 
 therefore it defines a point in the Siegel upper-half space

  $$\cH_g=\{ M\in GL(g,\bC) \ : \ {^t}\!M=M, \ \hbox{Im}(M)>0 \}.$$

The association $X \mapsto \Omega(X)$
can be interpreted as a {\em linearization} of the theory of compact 
Riemann surfaces  albeit an {\em imperfect} one
since, in order to compute the integrals, 
we made a choice of an integral basis (a lattice). 
The quotient of the Siegel space $\cH_g$ by this choice 
produces the moduli space of  
principally polarized complex abelian varieties 

$$
  \cA_g=\frac{\cH_g}{SP(2g,\bZ)}
$$

  and the association $X\mapsto \Omega(X)$ globalizes to a morphism (the {\em period map})

  $$ j:\cM_g\to \cA_g.$$
  
It is the right moment to mention two fundamental results 
of the last century Italian geometry school. The first one is the classic:
\begin{theorem}[Torelli, \cite{Torelli}] The map $j$ is injective.\end{theorem}
  
In other words, the period matrix reconstructs the curve;
 this should make the linearization effective but unfortunately $j$ is not surjective
 and, for $g>3$, the {\em Jacobian locus} $J_g=j(\cM_g)$ is not open in $\cA_g$.

The second one is the complete solution of the Schottky problem 
by Arbarello and De Concini in \cite{AdC}.  The famous Schottky problem consists in the
 characterization of the image $J_g\subset \cA_g$;
we lack the space for describing the solution in \cite{AdC} in details,
we only recall that it  
 has later been clarified by Shiota in \cite{S}. 

\subsection{ The form $\eta$}

Our starting point was the analysis of the local geometry of the inclusion $J_g\subset \cA_g$ 
of which the second fundamental form has been studied in \cite{CPT}.
The second fundamental form is a map

\begin{equation}\label{secondff}
\sff: Sym^2 T_{J_g,j([X])}\la N_{J_g,A_g,j([X])}
\end{equation}

where  $T_{J_g,j([X])}\equiv H^0(X,2K)^\ast$ is the tangent space 
to $J_g$ at a point $j([X])$ and $N$ is the normal space of the 
inclusion at $j([X])$. Recall that the dual of the normal corresponds 
to the quadrics containing the canonical curve. 

The idea, also supported by arithmetic conjectures 
(Coleman-Oort \cite{C} \cite{O} \cite{MO}), is that $J_g$ should be very curved.

In the work \cite{CPT}  an intrinsic map, named {\em Hodge-Gaussian}, has been constructed
and, at least in some particular cases, it allows explicit calculations.
\bigskip

A globalization of the Hodge-Gaussian map is constructed in 
\cite{CFG}.

On the product $X\times X$ we construct a meromorphic form $\eta$ 
with a double pole on the diagonal $\Delta$; more explicitly, if the diagonal has
 local equation $x=y$ then 

\begin{equation}\label{eta}
\eta(x,y)\equiv\frac {dx\wedge dy}{(x-y)^2}+ h(x,y)dx\wedge dy.
\end{equation}

with $h(x,y)$ holomorpic.

Some quick remarks are in order:
\begin{enumerate}
\item  $\eta$ is related to some sort of Green's function 
or, more precisely, to a Bergman kernel  (Ghigi and Tamborini, \cite{GT});

\item a different construction of $\eta$ appears in an unpublished book of Gunning \cite{Gu2};

\item many forms on $ X\times X$ with poles on the diagonal appear in the literature;
one of these even seems to be present in an ancient text of Klein.
 
 \end{enumerate}

All the information of the second fundamental form (\ref{secondff}) is captured by $\eta$. 
Indeed, we can rewrite all spaces involved in the construction of $\sff$  
by means of the variety $X\times X$; in particular there is a natural inclusion

 $$	\gamma: N_{J_g,A_g,j([X])}^\ast \hookrightarrow H^0 (X\times X,K_{X\times X}(-2\Delta)) $$ 
 
 of the dual of the normal into the space of holomorphic 
 $2$-forms with double zeroes on the diagonal.

\begin{theorem} [Colombo-Frediani-Ghigi, \cite{CFG}]
The restriction to $\gamma(N_{J_g,A_g,j([X])}^\ast)$ of the multiplication by the form $\eta$
$$ H^0(X\times X,K_{X\times X}(-2\Delta))  \stackrel{\cdot\eta}\la H^0(X\times X,2K_{X\times X})$$ 
is the dual of the second fundamental form $\sff$ induced by the period map.
 \end{theorem}

\subsection{Forms on the product and projective structures}

For any projective structure on $X$ we can construct a form on $X\times X$,
 defined in a neighbourhood of the diagonal $\Delta$ and 
 with a double pole on $\Delta$:  
start with a projective atlas $\{(U_\alpha,z_\alpha)\}_\alpha$ 
and set in $U_\alpha\times U_\alpha$, with coordinates
$ (x_\alpha, y_\alpha)$ induced by $z_\alpha$, the form 

\begin{equation} \frac{dx_\alpha\wedge dy_\alpha}{(x_\alpha-y_\alpha)^2}; \label{loc}\end{equation}
 this local expression gives a well-defined form on 
 $\bigcup_\alpha \left( U_\alpha \times U_\alpha\right) .$

\bigskip
Conversely, by solving a variant of the Schwarzian equation, 
we get a projective structure on $X$ from each of these forms 
(the story could be complex here too).

\begin{theorem}[Biswas-Raina \cite{BR}, Tyurin \cite{T}, Klein \cite{klein}] 

 Let  $W$ be a neighbourhood of the diagonal $\Delta\subset W\subset X\times X$, 
 and $\theta$ a differential form on $W$ with a double pole on $\Delta$
and biresidue $1$ (that is, it restricts to 1 on the diagonal), 
 $$\theta\in H^0(W,K_{X\times X}(2\Delta)),$$ 

 then there exists a unique projective structure $p(\theta)$ on $X$ such that,
  with respect to its projective atlas,  
  $\theta$ is expressed as in (\ref{loc}) up to second order.
\end{theorem}
At this point we can construct:
\begin{definition} The Hodge projective structure
$$p_h(X)$$ is the structure $p(\eta)$ induced by
 the differential form $\eta$ defined by $(\ref{eta}).$ 
\end{definition}
\bigskip


\section {The comparison between $p_u$ and $p_h$}

A natural problem is to decide whether the two canonical sections 
$p_u$ and $p_h$ are the same or not, 
and many suggested a positive response.

The result instead was:

\begin{theorem}[Biswas-Colombo-Frediani-Pirola (2021)] 
The two sections are not equivalent. That is, there exists an open subset $U\subset\cM_g$ 
such that for any $[X]\in U$ then $p_u(X)\neq p_h(X).$
\end{theorem}

The most interesting part is probably the method 
of the proof. In what follows, we briefly sketch it.

\bigskip
Recall that we have an exact sequence of vector spaces (see \cite{BR})
 $$0\to H^0(X,2K_X)\to E\stackrel{\phi} \to \bC\to 0$$ 
 and the projective structures on $X$ are the elements of
$$P(X)=\{ v\in E \ : \ \phi(v)=1\}.$$ 

By globalizing the sequence on $\cM_g$, we get 
an exact sequence of fiber bundles 

\begin{equation} 
0\to \Omega^1_{\cM_g}\to \cE\stackrel{\tilde\phi} \to \cO_{\cM_g}\to 0. \label{ex}
\end{equation} 
so that the sections corresponding to projective structures can be viewed as sections of $\cE$
 $$\cC^{\infty}(\cP)=\{\beta \in\cC^{\infty}(\cE) \ : \ \tilde\phi(\beta)=1\in H^0(\cM_g;\cO_{M_g}) \}.$$
 
 The exact sequence (\ref{ex}) does not split and the obstruction is a cohomology 
 class $\lambda$ which is a multiple of the generator of $H^2(\cM_g, \bZ).$ 
 This generator is the Hodge class.
Consider 
$$\overline{\partial} \beta $$
 and the computation 
$$\tilde\phi(\overline{\partial} \beta)=\overline{\partial} (1)=0.$$

Henceforth the form
$$\overline{\partial} \beta\in H^{1,1}(\cM_g)$$
 is a form of type $(1,1)$ that is, in local coordinates $z=(z_1, \dots, z_{3g-3})$
 can be expressed as 
 $$\overline{\partial} \beta=\sum a_{ij}(z)dz_i\wedge d\overline{z_j} $$

So, by construction:

\begin{proposition} 
With the previous notations we have that 
$$[\overline{\partial} \beta]=\lambda.$$
\end{proposition}

Moreover,

\begin{theorem} 
There is a $1$ to $1$ correspondence between the sections of the projective structures and
the $\overline{\partial}$-closed forms of type $(1,1)$ with Dolbeault cohomology class $\lambda$:
$$\cC^{\infty}(\cP)\longleftrightarrow \cA^{1,1}_{\lambda}(\cM_g)=\{\alpha\in \cA^{1,1}(\cM_g) \ : \  \overline{\partial} \alpha=0, \ [\alpha]=\lambda \}.$$
\end{theorem}

\begin{remark} The bijection above has been completed in a joint work 
with Favale and Torelli (\cite{FPT}) where it is shown that there are no 
holomorphic forms on $\cM_g.$ It should be recalled that Mumford showed 
in the sixties that there are no closed forms on $\cM_g$
\label{forme}\end{remark}

The comparison between $p_u$ and $p_h$ then becomes a comparison 
between the forms
 $$  \beta_u= \overline{\partial}(p_u) \quad \hbox{and} \quad \beta_h= \overline{\partial}(p_h).$$

Luckily the form $\beta_u$ is quite well-known:

\begin{theorem}\label{zttheorem}[Zograf-Takhtadzhyan, \cite{ZT}] 
The form $ \beta_u$ is the K\"{a}hler form of the 
Weil-Petersson metric: $$ \beta_u=\omega_{WP}.$$ 
In particular, it is everywhere positive definite.
\end{theorem}

The main point of our work was the calculation of $\beta_h$.

Recall the period map 
$$ j: \cM_g\to \cA_g $$
and the quotient 
$$\frac{\cH_g}{SP(2g,\bZ)}= \cA_g.$$

The Siegel space $\cH_g$ is a symmetric space and admits
 a natural metric with associated form $\omega_{Siegel}$ 
 which defines a form $\omega^\prime_{Siegel}$ on $\cA_g$.

\begin{theorem}[Biswas-Colombo-Frediani-Pirola, \cite{BCFP}] The form $ \beta_h$ 
is the pull-back of the Siegel metric via the period map:
$$\beta_h=j^\ast (\omega^\prime_{Siegel}) $$
\end{theorem}

\begin{proof}The (very simplified) idea beneath the proof 
is that both $ j^\ast (\omega^\prime_{Siegel})$ and $\beta_h$  are invariant 
under the action of the symplectic group $SP(2g,\bR)$; 
the uniqueness of the invariant forms concludes.
\end{proof}

The comparison between $j^\ast (\omega^\prime_{Siegel})$ 
and $\omega_{WP}$ is then easy:
when $g>2$, the form $j^\ast (\omega^\prime_{Siegel})$
is not positive definite in a sublocus where the differential $d j$ 
is not injective, namely the hyperelliptic locus; 
$\omega_{WP}$ is, on the contrary, everywhere positive definite.
In the genus $2$ case, an {\em ad hoc} argument is needed.

\section{Final remarks and open problems}

Recall the result in \cite{FPT} with Favale and Torelli where
 it is shown that there are no holomorphic forms on $\cM_g$
 and its corollary:

\begin {theorem} $H^0(\cM_g,\Omega^1)=0.$ \end{theorem}

\begin {corollary} $\cC^{\infty}(\cP)\equiv \cA^{1,1}_{\lambda}(\cM_g).$ \end{corollary}

Then, as a consequence (see \cite{BFPT}) we have 
a new point of view on the theorem of Zograf and Takhtadzhyan (\ref{zttheorem}).
 
There are three objects linked together by the study 
of the determinant of the Quillen cohomology:

 \begin{enumerate}
 \item the hyperbolic metric on a surface of genus $g>1$ (uniformization);
 \item the projective structure $p_u$;
 \item the Weil-Petersson metric.
\end{enumerate}

In \cite{BGT}, Biswas, Ghigi and Tamborini develop the same principle for $p_h$.
In this case the metric could be induced by the flat metric on the Jacobian of the curve (the Arakelov metric).

\bigskip

The geometric nature of the projective structure $p_h$ is still quite unclear,
as some natural questions are left unanswered:

 \begin{enumerate}
 \item determine the atlas, or at least understand a geometric construction;
 \item determine the monodromy and the associated fiber bundles: developing ({\em unrolling}) an atlas gives a monodromy action;
  \item determine the locus where $p_h=p_u$;
 \item what happens in some specific subloci of $\cM_g$?
\end{enumerate}

\section{Acknowledgments}

Both authors are very grateful to Alessandro Ghigi for his precoius remarks and some references.\\
 
The second author author is partially supported by PRIN project {\em Moduli spaces and special varieties} (2022). 

The first author acknowledges the projects `Fondo di Ateneo per la Ricerca 2019',  `Fondo di
Ateneo per la Ricerca 2020', funded by the University of Sassari. 

Both authors are members of GNSAGA, INdAM.

This work has been supported
by the National Research Project PRIN PNRR 20228CPHN5
“Metamaterials design and synthesis with applications to infrastructure engineering”
funded by the Italian Ministry of University and Research.

\end{document}